# The variational principle for weights characterizing the relevance

Mikhail A. Antonets, Grigoriy P. Kogan

1. **Introduction.**

The classical method of the thematic classification of texts is based on using the frequency weight on the list of words occurring in texts from the text corpus that determines the theme. In this method, the weight of each word is defined as its normalized frequency in the corpus's texts. The frequency weight is applied for determining the relevance of the tested text to the theme given via a text corpus: the tested text's relevance is defined as the value of its frequency weight (see [1]-[3]). In the present work we propose a method of constructing some optimal weights generated by certain variational principles leading to LP (linear programming) problems. A noteworthy feature of those optimal weights is a relatively small number of words belonging to their supports: in many the examples we considered that number didn't exceed 10% of the quantity of different words (lemmas) occurring in the corpus's texts, and besides the majority of those words were apparent candidates to be selected as the corpus's key words. The application of our method to the determination of the relevance of texts to a given collection of thematic text corpuses demonstrated its high efficiency and performance.

In the present paper, we offer a new method of constructing weights (probability measures) on a finite set $V$ that characterize the relevance of an arbitrary non-negative function on $V$ to a given set $M$ of non-negative functions on $V$.

As a quantitative estimation of the degree of relevance, i.e. the correspondence of the tested function $f$ to the criteria for the formation of the set of functions $M$ is based on, we propose the value

$$\underline{r}^M(f) = \sum_{v \in V} f(v) \underline{x}^M(v) \qquad (1)$$

where $\underline{x}^M$ is the constructed weight and

$$|f| = \sum_{v \in V} f(v) \quad .$$

Another weight $\overline{x}^M$ characterizes the irrelevance of the function $f$ by means of the value

$$\overline{s}^M(f) = \sum_{v \in V} f(v) \overline{x}^M(v) \qquad (2)$$

which decreases with the relevance's growth.

The method of constructing weights is based on certain variational principles leading to the problems of minimizing the corresponding linear functionals on a polyhedron.

The proposed formulation of the problem of estimating the relevance can be applied in the analysis of experimental data. As an example of such an application, let's consider the problem of classifying complex events assumed to be cycles of controlled technical systems' work. In this example, the set $V$ is the set of control actions accomplished in the process of the system management, while every function $m$ from the set $M$ takes a nonnegative value $m(v)$ on each of the control actions. The numbers $m(v)$ may be the values of controlling parameters, the switch numbers of 2-positioned relays, or other values characterizing the management process.

If the set $M$ consists of functions received in the studied managed system's cycles that were recognized as acceptable by experts then the proposed method allows to find a definite answer to the following question:

Whether the work cycle corresponding to the given function $f$ on the set $V$ is acceptable?

## 2. Mathematical background.

**Definition 1.** A weight on a finite set $V$ is a non-negative function $x \in \mathbb{R}^V$ satisfying the equality

$$|x| = 1 \qquad (3)$$

The set of all weights on the set $V$ will be denoted by $\Delta_V$.

Let $M$ be a finite set of nonnegative functions from the space $\mathbb{R}^V$.
For any $x, m \in \mathbb{R}^V$ we denote

$$(x,m)_V = \sum_{v \in V} x(v) m(v)$$

**Definition 2.** A weight $\underline{x}^M$ on a set $V$ will be called **supporting** for a set $M$ of nonnegative functions on the set $V$ if for any weight $x$ on the set $V$ it satisfies the inequality

$$min_{m \in M} (\underline{x}^M, m)_V \geq min_{m \in M} (x, m)_V$$

The set of supporting weights $\{\underline{x}^M\}$ coincides with the set of solutions of the variational problem

$$min_{m \in M} (\underline{x}^M, m)_V = max_{x \in \Delta_V} min_{m \in M} (x, m)_V \qquad (4)$$

**Definition 2.** A weight $\overline{x}^M$ on a set $V$ will be called **covering** for a set $M$ of nonnegative functions on the set $V$ if for any weight $x$ on the set $V$ it satisfies the inequality

$$max_{m \in M} (\overline{x}^M, m)_V \leq max_{m \in M} (x, m)_V$$

The set of covering weights $\{\overline{x}^M\}$ coincides with the set of solutions of the variational problem

$$max_{m \in M} (\overline{x}^M, m)_V = min_{x \in \Delta_V} max_{m \in M} (x, m)_V \qquad (5)$$

The supporting and covering weights are well known *maxmin* and *minimax* probability distributions respectively (see [4] for instance).

The below-formulated Theorem 1 states that the both variational problem can be reduced to the following well-known linear programming problems:

**Problem A.** Find all the pairs $\{\underline{\alpha}^M, \underline{x}^M\}$, $\underline{\alpha}^M \in \mathbb{R}$, $\underline{x}^M \in \Delta_V$, that maximize the value $\alpha$ under the conditions (3), inequalities

$$x(v) \geq 0, v \in V ,  \qquad (6)$$

and inequalities

$$\alpha - (x, m)_V \leq 0, \quad \forall m \in M \qquad (7)$$

**Problem B.** Find all the pairs $\{\overline{\alpha}^M, \overline{x}^M\}$, $\overline{\alpha}^M \in \mathbb{R}$, $\overline{x}^M \in \Delta_V$, that maximize the value $\alpha$ under the condition, inequalities

$$x(v) \geq 0, v \in V ,  \qquad (8)$$

and inequalities

$$\alpha - (x, m)_V \geq 0, \quad \forall m \in M \qquad (9)$$

For an arbitrary element $v \in V$, let's denote by $\hat{v}$ the function on the set $M$ defined through the relationship $\hat{v}(m) = m(v)$ for any $m \in M$ and let's also denote $\hat{V} = \{\hat{v},\ v \in V\}$.

Therefore we can construct the supporting weight $\underline{x}^{\hat{V}}$ and the covering weight $\overline{x}^{\hat{V}}$ for the set $\hat{V}$ of functions on the set $M$. These weights can be found as solutions of problems $\hat{A}, \hat{B}$ presented below.

**Problem $\hat{A}$.** Find all the pairs $\{\underline{\alpha}^{\hat{V}}, \underline{x}^{\hat{V}}\}$, $\underline{\alpha}^{\hat{V}} \in \mathbb{R}$, $\underline{x}^{\hat{V}} \in \Delta_M$, that maximize the value $\alpha$ for the condition of fulfilling the equality (3), inequalities

$$x(m) \geq 0, m \in M ,  \qquad (10)$$

and inequalities

$$\alpha - (x, \hat{v})_M \leq 0, \quad \forall v \in V \qquad (11)$$

**Problem $\hat{B}$.** Find all the pairs $\{\overline{\alpha}^{\hat{V}}, \overline{x}^{\hat{V}}\}$, $\overline{\alpha}^{\hat{V}} \in \mathbb{R}$, $\overline{x}^{\hat{V}} \in \Delta_M$, that maximize the value $\alpha$ for the condition of fulfilling the equality (3), inequalities

$$x(m) \geq 0, m \in M ,  \qquad (12)$$

and inequalities

$$\alpha - (x, \hat{v})_M \geq 0, \quad \forall v \in V \qquad (13)$$

**Theorem 3.** The sets of supporting weights and covering weights for the sets of functions $M$ and $\hat{V}$ are non-empty and there exist the following connection between them:

1. Linear programming problem $\hat{B}$ is dual problem for linear problem $A$ and therefore
$$\overline{\alpha}^{\hat{V}} = \underline{\alpha}^M$$
2. Linear programming problem $\hat{A}$ is dual problem for linear problem $B$ and therefore

$$\overline{\alpha}^M = \underline{\alpha}^V$$

3. The complementary slackness condition gives the relations

$$\sum_{v \in V} m(v)\underline{x}^M(v) = \underline{\alpha}^M \quad \text{if} \quad \overline{x}^V(m) \neq 0 \qquad (14)$$

$$\sum_{m \in M} m(v)\overline{x}^V(m) = \underline{\alpha}^M \quad \text{if} \quad \underline{x}^M(v) \neq 0 \qquad (15)$$

$$\sum_{v \in V} m(v)\overline{x}^M(v) = \underline{\alpha}^V \quad \text{if} \quad \underline{x}^V(m) \neq 0 \qquad (16)$$

$$\sum_{m \in M} m(v)\underline{x}^V(m) = \underline{\alpha}^V \quad \text{if} \quad \overline{x}^M(v) \neq 0 \qquad (17)$$

4. From the above relationships, there follow that

$$\underline{r}^M(m) \geq \underline{\alpha}^M \qquad (18)$$

$$\overline{s}^M(m) \leq \overline{\alpha}^M \qquad (19)$$

and any function $f$ on the set $V$ is relevant to the set $M$ if

$$\underline{r}^M(f) \geq \underline{\alpha}^M$$

### 3. Examples.

Now let's consider some examples of supporting and covering weights.

**Example 1.** Let the function $m$ from the set $M$ be defined by the following relationship

$$m(v) = \begin{cases} 1, & v = v_m \\ \varepsilon, & v \neq v_m \end{cases}$$

where $\varepsilon > 0$ and events $v_m$ are different for different $m$, that is possible when $|V| \geq |M|$. Let's denote by $V_M$ the set of elements from $V$ such that in each of them one of the functions from the set $M$ equals 1.

In the latter case the supporting weight $\underline{x}^M$ is a solution of the following system of the inequalities:

$$x(v_m) + \varepsilon(1 - x(v_m)) \geq \alpha, \forall v_m \in V_M$$

for biggest possible value of the variable $\alpha$.

The symmetry of this system of inequalities implies that

$$\underline{x}^M(v) = \frac{1}{|V_M|}, \forall v \in V_M$$

and

$$\underline{\alpha}^M = \frac{1}{|V_M|} + \varepsilon\left(1 - \frac{1}{|V_M|}\right)$$

Acting by analogy, in the covering weight problem we receive the equalities

$$\overline{\alpha}^M = \underline{\alpha}^M$$

$$\overline{x}^M(v) = \underline{x}^M(v), \forall v \in V_M$$

Hence, in the last case the relevance measure coincides with the irrelevance one.

For the first view, it seems paradoxical. But, in fact, it means that the functions from the set $M$ doesn't have common structural singularities and, therefore, the question about relevance is senseless in such a case.

**Example 2.** Let all the functions from the set $M$ be equal to 1 in the points of their support and some point $v_0$ be the only common support point for any pairs of functions from $M$, and besides for each function $m \in M$ let the sets $V_m = supp\ m \setminus \{v_0\}$ be non-empty.

Then the relationships (6) for the set $M$ obtain the form

$$\alpha - x(v_0) - \sum_{v \in V_m} x(v) \leq 0, \quad \forall m \in M$$

and the only non-negative solution of these inequalities having the maximal value of the variable $\alpha$ has the form

$$\underline{\alpha}^M = 1,$$

$$\underline{x}^M(v_0) = 1,\ \underline{x}^M(v) = 0\ \text{для}\ v \neq v_0$$

In the considered case, the inequalities determining the covering weights $\overline{x}^M$ have the form

$$\alpha - x(v_0) - \sum_{v \in V_m} x(v) \geq 0$$

and they imply that for the minimal $\alpha$

$$x(v_0) = 0$$

and

$$\overline{\alpha}^M = \sum_{v \in V_m} \overline{x}^M(v) = \frac{1}{|M|}$$

From the latter relationship, we can conclude that there exists a continuum of covering weights if at least one of the sets $V_m$ contains more than one element.

### References.